\documentclass[11pt]{article}
\usepackage{amsmath,amsthm,amsfonts,amssymb,amscd, amsxtra}

\newcommand{\banacha}{X}
\newcommand{\banachb}{Y}
\newtheorem{theorem}{Theorem}
\newtheorem{lemma}[theorem]{Lemma}

\newtheorem{corollary}[theorem]{Corollary}
\newtheorem{proposition}[theorem]{Proposition}

\newtheorem{remark}{Remark}
\newtheorem{example}{Example}
\begin{document}
\title{Local convergence analysis of  inexact Newton-like \\  methods under  majorant condition\thanks{IME/UFG, Campus II- Caixa
    Postal 131, CEP 74001-970 - Goi\^ania, GO, Brazil.}}

\author{ O. P. Ferreira\thanks{E-mail:{\tt
      orizon@mat.ufg.br}. The author was supported in part by
     CNPq Grant 302618/2005-8,  PRONEX--Optimization(FAPERJ/CNPq) and FUNAPE/UFG.}
  \and M. L. N. Gon\c calves \thanks{E-mail:{\tt
      maxlng@hotmail.com}.The author was supported in part by CAPES.}
}
\date{July 23, 2008}

\maketitle
\begin{abstract}
We present a local convergence analysis of inexact Newton-like methods for solving nonlinear equations under majorant conditions. This analysis provides an estimate of the convergence radius and a clear relationship between the majorant function, which relaxes the Lipschitz continuity of the derivative, and the nonlinear operator under consideration.  It also allow us to obtain  some important special cases.\\


\noindent
{\textsc AMSC: 49M15, 90C30.}
\end{abstract}
\section{Introduction}\label{sec:int}
Newton's method and its variations (see \cite{DN1}), including the inexact Newton methods, are the most efficient methods known for solving nonlinear equations
\begin{equation} \label{eq:nle}
F(x)=0,
\end{equation}
where  $F:{\Omega}\to \mathbb{R}^{n}$ is a continuously differentiable function  and $\Omega\subseteq \mathbb{R}^{n}$ is an open set. The inexact Newton method was introduced by Dembo, Eisenstat and Steihaug in \cite{DE1} denoting any method which, given an initial point $x_0$, generates the sequence $\{ x_k\}$  as follows:
$$
x_{k+1}={x_k}+S_k,\qquad k=0,1,\ldots,
$$
where $S_k$ is the solution of the linear system
$$
F'(x_k)S_k=-F(x_k)+r_{k},
$$
for a suitable residual  $r_k\in \mathbb{R}^{n}$. Let $x_{*}$  be a solution of \eqref{eq:nle} such that $F'(x_*)$ is invertible. As shown in \cite{DE1}, if $\|r_{k}\|\leq \theta_{k}\|F(x_{k})\|$ for $k=0,1,\ldots$ and $\{\theta_{k}\}$ is a sequence of forcing terms such that $0\leq \theta_{k}<1$ then there exists $\epsilon>0$ such that the sequence $\{ x_k\}$, for any initial point $x_0 \in B(x_*, \epsilon)=\{x\in \mathbb{R}^{n}:\; \|x_{*}-x\|<\epsilon\}$, is well defined and converges linearly to  $x_{*}$  in the  norm $\|y\|_*=\|F'(x_*)y\|$, where $\| \; \|$ is any norm in $ \mathbb{R}^{n}$. It is worth noting that, in \cite{DE1}, no Lipschitz condition is assumed on the derivative $F'$ to prove that $\{x_{k}\}$ is well defined and linearly converging; however, no estimate of the convergence radius $\epsilon$ is provided. As pointed out by \cite{JM10} (see also \cite{B10})   the result of \cite{DE1} is difficult to apply due to dependence of the norm  $\|\; \|_*$, which is not computable.

Formally, the inexact Newton-like methods for solving the non-linear equation \eqref{eq:nle}, which we will consider,  are described as follows: Given an initial point $x_0 \in {\Omega}$, define
$$
x_{k+1}={x_k}+S_k,\qquad B(x_k)S_k=-F(x_k)+r_{k}, \qquad k=0,1,\ldots,
$$
where $B(x_k)$ is a suitable invertible approximation of the derivative $F'(x_k)$ and the  residual  $r_k$ satisfies 
$$
\|P_{k}r_{k}\|\leq \theta_{k}\|P_{k}F(x_{k})\|,
$$
for suitable forcing sequence $\{\theta_{k}\}$ and some invertible matrix sequence $\{P_{k}\}$ of preconditioners for the above linear equation defining the steep $S_k$. This method was considered for the first time in  \cite{B10}, and was also considered in \cite{C10}. In particular, letting $P_k\equiv I$ be the identity matrix  and $B_k=F'(x_k)$ for each $k$, we obtain the inexact Newton method considered in \cite{DE1}, \cite{Mo89} and \cite{Y84}.

Inexact Newton-like methods may fail to converge and may even fail to be well defined. To ensure that the method is well defined and converges to a solution of a given non-linear equation, some conditions must be imposed. For instance, the classical convergence analysis (see \cite{DE1}) requires the initial iterate to be "close enough" to a solution and the first derivative of the non-linear function to be invertible in this solution. Moreover, for estimating the convergence radius, the Lipschitz continuity or something like Lipschitz continuity, of the first derivative is also assumed (see   \cite{C10}, \cite{B10} and \cite{Y84}).

In the last twenty-five years, there have been papers dealing with the issue of convergence of the Newton methods, including the inexact Newton methods and the Gauss-Newton methods, by relaxing the assumption of Lipschitz continuity of the derivative (see  \cite{ABM2004}, \cite{F08}, \cite{FS06},  \cite{C08}, \cite{C10}, \cite{Mo89},  \cite{XW10}, \cite{Wu08} and \cite{Y84}). In addition to improving the convergence theory (this allows us to estimate the convergence radius and to enlarge the range of application) some modifications of the Lipschitz condition also permit us to unify several results. Works dealing with this subject include \cite{ABM2004}, \cite{F08}, \cite{FS06} and \cite{Wu08}.


Our aim in this paper is to present a new local convergence analysis for inexact Newton-like methods under majorant condition. In our analysis, the classical Lipschitz condition is relaxed using a majorant function. It is worth pointing out that this condition is equivalent to Wang's condition introduced in \cite{XW10} and used by Chen and Li in \cite{C10} to study the inexact Newton-like methods. The convergence analysis presented is linear in an arbitrary norm. It provides a new estimate for the convergence radius  and  a clear relationship between the majorant function and the nonlinear operator under consideration. It also allows us to obtain some special cases that can be evaluated as an application.

The organization of the paper is as follows.
In Section \ref{sec:int.1}, we list some notations and basic results used in our presentation.
In Section \ref{lkant} the main result is stated, and in  Section \ref{sub:mf} some properties involving the majorant function are established. In Section \ref{sec:MFNLO} we presented the relationships between the majorant function and the non-linear operator. In Section \ref{sec:proofqn} the main result is proveda and  some applications of this result are given in Section \ref{sec:ec}. Some final remarks are made in Section~\ref{rf}.

\subsection{Notations and  auxiliary  results} \label{sec:int.1}
The following notations and results are used throughout our presentation. Let $\mathbb{R}^n$ be with a norm $\| . \|$. The open and closed ball
at $a \in \mathbb{R}^n$ and radius $\delta>0$ are denoted, respectively by
$$ 
B(a,\delta) =\{ x\in \mathbb{R}^n ;\; \|x-a\|<\delta \}, \qquad B[a,\delta] =\{ x\in \mathbb{R}^n ;\; \|x-a\|\leqslant \delta \}.
$$
Let $\mathcal{L}(\mathbb{R}^n, \mathbb{R}^n)$ be the space of  liner operators of $\mathbb{R}^n$. Define the operator norm associated to the norm  $\| . \|$ as 
$$
\|T\|:=\sup\{ \|Tx\|,\;\|x\|\leq 1\},  \qquad T\in \mathcal{L}(\mathbb{R}^n, \mathbb{R}^n).
$$
The condition number of an invertible operator  $T$  is denoted by $\mbox{cond}(T):=\|T\|\|T^{-1}\|.$
\begin{lemma}[Banach's Lemma] \label{lem:ban}
Let $B \in \mathcal{L}(\mathbb{R}^n, \mathbb{R}^n)$   and   $I\in \mathcal{L}(\mathbb{R}^n, \mathbb{R}^n)$ , the identity
operator. If  $\|B-I\|<1$,  then $B$
is invertible and  $ \|B^{-1}\|\leq 1/\left(1-
\|B-I\|\right). $
\end{lemma}
\begin{proof}  See the proof of Lemma 1, pp. 189 of Smale \cite{S86}  with  $A=I$  and  $c=\|B-I\|$.
\end{proof}

\begin{proposition} \label{le:ess}
If $0\leq t <1$, then $ \sum
_{i=0}^{\infty}(i+2)(i+1)t^{i}=2/(1-t)^3. $
\end{proposition}
\begin{proof} Take $k=2$ in  Lemma 3, pp. 161 of  Blum, Cucker,  Shub and Smale \cite{BCSS97}.
\end{proof}
Also, the following auxiliary result of elementary convex analysis
will be needed:
\begin{proposition}\label{pr:conv.aux1}
Let $\epsilon  >0$  and  $\tau \in [0,1]$. If $\varphi:[0,\epsilon) \rightarrow\mathbb{R}$ is convex, then $l:(0,\epsilon) \to \mathbb{R}$ defined by 
$$
l(t)=\frac{\varphi(t)-\varphi(\tau t)}{t},
$$
is non-increasing.  
\end{proposition}
\begin{proof} See  Theorem 4.1.1 and Remark 4.1.2 on pp. 21 of Hiriart-Urruty and Lemar\'echal \cite{HL93}.
\end{proof}
\section{ Local convergence of inexact Newton-like method} \label{lkant}

Our goal is to state and prove a local theorem for {\it inexact Newton-like methods}. Assuming that the non-linear equation $F(x)=0$ has a solution $x_*$, we will, under mild conditions, prove that the inexact Newton-like method is well defined and that the generated sequence converges linearly to this solution. The statement of the theorem is as follows:
\begin{theorem}\label{th:ntqn}
Let $\Omega\subseteq \mathbb{R}^{n}$ be an open set and $F:{\Omega}\to \mathbb{R}^{n}$ a continuously differentiable function. Let  $x_*\in \Omega$,  $R>~0$ and 
$$
\kappa:=\sup \left\{ t\in [0, R): B(x_*, t)\subset\Omega \right\}.
$$
Suppose that $F(x_*)=0$, $F '(x_*)$ is invertible and there exists a $f:[0,\; R)\to \mathbb{R}$ continuously differentiable such that 
  \begin{equation}\label{Hyp:MHqn}
\left\|F'(x_*)^{-1}\left[F'(x)-F'(x_*+\tau(x-x_*))\right]\right\|
\leq
    f'\left(\|x-x_*\|\right)-f'\left(\tau\|x-x_*\|\right),
  \end{equation}
  for $\tau \in [0,1]$ and  $x\in B(x_*, \kappa)$, where
\begin{itemize}
  \item[{\bf h1)}]  $f(0)=0$ and $f'(0)=-1$;
  \item[{\bf  h2)}]  $f'$ is  convex and strictly increasing.
\end{itemize}
Take $0\leq \vartheta<1$, $0\leq \omega_{2}<\omega_{1}$  such that  $\omega_{1}\vartheta+\omega_{2}< 1$.  Let
 $
 \nu:=\sup\{t\in [0, R): f'(t)<0\},
 $
 $$\rho:=\sup\{t\in(0, \nu):\omega_{1}(1+\vartheta)[f(t)/(tf'(t))-1]+\omega_{1}\vartheta+\omega_{2}<1\}, \qquad
 \sigma:=\min \left\{\kappa, \, \rho\right\}.
$$
Then, the inexact Newton-like methods for solving $F(x)=0$, with initial point $x_0\in
B(x_*, \sigma)\backslash \{x_*\}$
\begin{equation} \label{eq:DNSqn}
x_{k+1}={x_k}+S_k, \qquad  B(x_k)S_k=-F(x_k)+r_{k}, \qquad
\; k=0,1,\ldots,
\end{equation}
where $B(x_k)$ is an invertible approximation of $F'(x_k)$ satisfying 
\begin{equation*}\label{con:qn}
\|B(x_k)^{-1}F'(x_k)\| \leq \omega_{1}, \qquad
\|B(x_k)^{-1}F'(x_k)-I\| \leq \omega_{2}, 
\end{equation*}
the  residual  $r_k$ satisfies 
\begin{equation}\label{eq:ERROqn}
\|P_{k}r_{k}\|\leq \theta_{k}\|P_{k}F(x_{k})\|,
\end{equation}
for some invertible matrix sequence $\{P_{k}\}$ of preconditioners (for the linear system in \eqref{eq:DNSqn}) and a forcing sequence $\{\theta_{k}\}$ of non-negative numbers satisfying 
\begin{equation*}\label{con:vk}
\theta_{k}\mbox{cond}(P_{k}F'(x_{k}))\leq
\vartheta,
\end{equation*}
is well defined, contained in $B(x_*,\sigma)$, converges to $x_*$ and there holds
  \begin{equation*} \label{eq:q2qn}
    \|x_{k+1}-x_*\| \leq
    \left[\omega_{1}(1+\vartheta)\left(\frac{f(\|x_{0}-x_*\|)}{\|x_{0}-x_*\|
    f'(\|x_{0}-x_*\|)}-1\right)+
    \omega_{1}\vartheta+\omega_{2}\right]\|x_k-x_*\|, \quad k=0,1,\ldots .
  \end{equation*}
 \end{theorem}
 Note that letting  the majorant function $f:[0, \kappa)\to \mathbb{R}$ be given by 
$
f(t)=Kt^{2}/2-t,
$
$B({x_k})=F'(x_k)$, $P_k=I$,   $\omega_1=1$ and
$\omega_2=0$ in  Theorem \ref{th:ntqn}, we obtain the following local convergence result for the inexact Newton method:
 \begin{theorem}\label{th:pclc}
Let $\Omega\subseteq \mathbb{R}^{n}$ be an open set and
$F:{\Omega}\to \mathbb{R}^{n}$ be continuously differentiable in
$\Omega$. Take  $x_*\in \Omega$ and let
$
\kappa:=\sup \left\{ t>0: B(x_*, t)\subset\Omega \right\}.
$
Assume that $F '(x_*)$ is invertible, $F(x_*)=0$, there exists a
$K>0$  such that
  $$
\left\|F'(x_*)^{-1}\left[F'(x)-F'(y)\right]\right\| \leq K\|x-y\|, \qquad \forall\; x, y\in B(x_*, \kappa).
$$
Take $0\leq \vartheta<1$. Let
$$\sigma:=\min\left\{\kappa,
\,2(1-\vartheta)/\left(K(3-\vartheta)\right)
\right\}.$$
 Then, the inexact Newton-like method for solving
$F(x)=0$, with the initial point $x_0\in B(x_*, \sigma)\backslash \{x_*\}$
$$
x_{k+1}={x_k}+S_k, \qquad  F'(x_k)S_k=-F(x_k)+r_{k}, \qquad
\; k=0,1,\ldots,
$$
where the  residual  $r_k$ satisfies
$$
\|r_{k}\|\leq \theta_{k}\|F(x_{k})\|,
$$
for some forcing sequence $\{\theta_{k}\}$   of non-negatives numbers, satisfying 
$$
\theta_{k}\mbox{cond} F'(x_{k})\leq \vartheta,
$$
is well defined, contained in $B(x_*,\sigma)$, converges to $x_*$
and there holds
$$
    \|x_{k+1-x_*}\| \leq
\left[(1+\vartheta)\frac{K\|x_{0}-x_*\|}{2(1-K\|x_{0}-x_*\|)}+\vartheta
\right]\|x_k-x_*\|, \quad \forall \; k=0,1,\ldots .
 $$
 \end{theorem}
 \begin{remark}
Note that letting   $\vartheta=0$ (in this case $\theta_k\equiv 0$ and $r_k\equiv 0$) in Theorem~\ref{th:pclc} we  obtain $r=\min \left\{\kappa, \,2/(3K) \right\}$. 
As was shown in   \cite{TW79} (see also, \cite{Y82}) this is the best possible convergence radius for Newton's Method. Therefore, for vanishing residuals, Theorem~\ref{th:pclc} merges into the theory of Newton's Method and, as a consequence,  Theorem~ \ref{th:ntqn} does too.
\end{remark}

In order to prove Theorem \ref{th:ntqn} we need some results. From here on, we assume  that all assumptions of Theorem \ref{th:ntqn} hold.
\subsection{The majorant function} \label{sub:mf}
Our first goal is to show that the constant $\kappa$ associated with $\Omega$ and the constants  $\nu$,  $\rho$ and $\sigma$ associated with the majorant function $f$ are positive. Also, we will prove some results related to the function $f$.

We begin by noting that  $\kappa>0$, because $\Omega$ is an open set and $x_*\in \Omega$.
\begin{proposition} \label{pr:incr1}
The following statements hold:

\begin{itemize}
 \item[{\bf i)}] $\nu>0$;
 \item[{\bf ii)}] $f'(t)<0, \quad \forall \; t\in [0,\,\nu);$
\item[{\bf iii)}] The map $[0,\, \nu) \ni t \mapsto 1/|f'(t)|$
is strictly  increasing and
 \item[{\bf iv)}] $ t-f(t)/f'(t)<0,\quad \forall \; t\in [0,\,\nu).$
\end{itemize}
\end{proposition}
\begin{proof} 
 As $f'$ is continuous in $(0, R)$ and $f'(0)=-1$, there exists a $\delta>0$ such that $f'(t)<0$ for all $t\in (0,\, \delta).$ So,  $\nu\geq \delta$ and item~{\bf i} is proved. 
 
For proving item~{\bf ii}, use {\bf h2} and the definition of $\nu$. Now, for establishing the validity of item~{\bf iii}, combine {\bf h2} and item~{\bf ii}.
 
Since  $f'$ is strictly increasing we have  $f$ is strictly convex. So, 
$$
f(0)>f(t)-tf'(t), \qquad \forall \;t\in  (0,\, R).
$$
Because $f(0)=0$ and $f'(t)<0$ for all $t\in (0, \,\nu)$, the inequality in item~{\bf iv}  follows from above equation.
\end{proof}

Let $n_{f}$ be the Newton iteration associated with the majorant function,  
\begin{equation} \label{eq:def.nf}
  \begin{array}{rcl}
  n_{f}:[0,\, \nu)&\to& (-\infty, \, 0]\\
    t&\mapsto& t-f(t)/f'(t).
  \end{array}
\end{equation}
From Proposition~\ref{pr:incr1} we have   $f'< 0$ in $[0, \,\nu)$.
Hence, the Newton iteration associated with the majorant function is well defined in $[0,\, \nu)$.
\begin{proposition}  \label{pr:incr2}
The map  $(0,\, \nu) \ni t \mapsto |n_{f}(t)|/t^2$  is  strictly  increasing.
\end{proposition}
\begin{proof}
Using item~{\bf iv} of Proposition~\ref{pr:incr1} and {\bf h1} we obtain, after simple algebraic manipulation,  that 
\begin{equation} \label{eq:deta}
 \frac{|n_{f}(t)|}{t^2}=
\frac{1}{|f'(t)|}\int_{0}^{1}\frac{f'(t)-f'(\tau t)}{t} \,d
\tau, \qquad \forall \;t\in (0,\, \nu).
\end{equation}
On the other hand, since $f'$ is strictly increasing, we obtain that the map
$$
(0,\, \nu) \ni t \mapsto \frac{f'(t)-f'(\tau t)}{t},
$$
is positive for all $\tau
\in (0, 1)$. Also, from {\bf  h2} we know that  $f'$ is convex. So, applying 
Proposition~\ref{pr:conv.aux1} with $f'=\varphi$ and
$\epsilon=\nu$, we conclude that the last map is increasing. Hence the second term in the right hand side of \eqref{eq:deta} is positive and increasing. Therefore, since Proposition~\ref{pr:incr1} implies that the first term in the right had side of \eqref{eq:deta} is positive and strictly  increasing, we conclude the statement.

\end{proof}
\begin{corollary}  \label{pr:10}
The map $(0,\, \nu) \ni t \mapsto |n_{f}(t)|/t$  is  strictly  increasing.
\end{corollary}
\begin{proof}It is immediate, by noting that 
$|n_{f}(t)|/t=(|n_{f}(t)|/t^2)t$ is a product of two strictly  increasing functions.
\end{proof}
\begin{proposition}\label{pr:incr102}
The constant $ \rho $ is positive and there holds
$$
\omega_{1}(1+\vartheta)\frac{|n_{f}(t)|}{t}+
\omega_{1}\vartheta+\omega_{2}<1, \qquad
\forall \; t\in (0, \, \rho).
$$
\end{proposition}
\begin{proof}
Using Proposition \ref{pr:incr1} and the definition  \eqref{eq:def.nf},
we have
\begin{equation} \label{eq:rho1qn}
0<f(t)/(tf'(t))-1=\big(f(t)/f'(t)-t\big)/t=|n_{f}(t)|/t,
\qquad \forall \; t\in (0,\, \nu).
\end{equation}
Now, because Proposition \ref{pr:incr2} implies
that $|n_{f_r}(t)|/t^2$ is bounded near zero, we obtain 
\begin{equation}\label{t12}
\lim_{t\to 0}|n_{f}(t)|/t=\lim_{t\to 0}(|n_{f}(t)|/t^2)\,t=0.
\end{equation}
  Thus, since $
1-(\omega_{1}\vartheta+\omega_{2})/\omega_{1}(1+\vartheta)>0, $
using \eqref{eq:rho1qn} and  \eqref{t12} we conclude that there exists a
$\delta>0$ such that
$$
0<(f(t)/(tf'(t))-1)<1-(\omega_{1}\vartheta+\omega_{2})/\omega_{1}(1+\vartheta),
 \qquad \forall \; t\in (0, \delta),
$$
or, equivalently, 
\begin{equation} \label{eq:rho2qn}
0<\omega_{1}(1+\vartheta)[f(t)/(tf'(t))-1]+\omega_{1}\vartheta+\omega_{2}<1,
\qquad \forall \; t\in (0, \delta).
\end{equation}
Hence, combining the last equation and the definition of $\rho$, we have  $\delta\leq \rho$, which is a  proof of the first statement.

For concluding the proof, we use the definition of $\rho$,  equality \eqref{eq:rho2qn},  \eqref{eq:rho1qn} and Corollary \ref{pr:10}.
\end{proof}
\subsection{Relationship of the majorant function with the non-linear operator} \label{sec:MFNLO}
In this section we will present the main relationships between the  majorant function $f$ and the non-linear operator $F$. 
\begin{lemma} \label{wdns}
Let $x \in \Omega$. If \,\,$\| x-x_*\|<\min\{\nu, \kappa\}$,  then $F'(x) $ is invertible and
$$
\|F'(x)^{-1}F'(x_*)\|\leqslant  1/|f'(\| x-x_*\|)|. 
$$
In particular, $F'$ is invertible in $B(x_*, \sigma)$.
\end{lemma}
\begin{proof}
Let $x\in \Omega$ such that $\| x-x_*\|<\min\{\nu, \kappa\}$. So $f'(\|x-x_*\|)<0$ which, together with \eqref{Hyp:MHqn}, implies 
  \begin{align*}
    \|F'(x_*)^{-1}F'(x)-I\|= \|F'(x_*)^{-1}[F'(x)-F'(x_*)]\|
    &\leq f'(\|x-x_*\|)-f'(0)<-f'(0)=1.
  \end{align*}
Thus, Lemma \ref{lem:ban} and the last equation imply that
  $F'(x_*)^{-1}F'(x)$ is invertible, as well as $F'(x)$,  and
  $$
    \|F'(x)^{-1}F'(x_*)\| \leq
    \frac{1}{1-\|F'(x_*)^{-1}F'(x)-I\|}\leq
    \frac{1}{1-\left(f'(\| x-x_*\|)-f'(0)\right)}= \frac{1}{|f'(\| x-x_*\|)|},
  $$
where we assume that $f'(0)=-1$ and $f'<0$ in $[0, \, \nu)$ in the last equality. As $\sigma\leq \nu$ the last part is proved.
\end{proof}
The Newton iteration at a point happens to be a zero of the linearization
of $F$ at such point, which is also the first-order Taylor expansion
of $F$.  So, we study the linearization error at point
in $\Omega$
\begin{equation}\label{eq:def.er}
  E_F(x,y):= F(y)-\left[ F(x)+F'(x)(y-x)\right],\qquad y,\, x\in \Omega.
\end{equation}
We will bound this error by the error of the linearization of the
majorant function $f$
\begin{equation}\label{eq:def.erf}
        e_f(t,u):= f(u)-\left[ f(t)+f'(t)(u-t)\right],\qquad t,\,u \in [0,R).
\end{equation}
\begin{lemma}  \label{pr:taylor}
If  $\|x_*-x\|< \kappa$, then there holds
$
\|F'(x_*)^{-1}E_F(x, x_*)\|\leq e_f(\|x-x_*\|, 0). 
$
\end{lemma}
\begin{proof}
 Since   $B(x_*, \kappa)$ is convex, we obtain that $x_*+\tau(x-x_*)\in B(x_*, \kappa)$, for $0\leq \tau \leq 1$. Thus, as $F$ is continuously differentiable in $\Omega$, the definition of $E_F$ and some simple manipulations yield
$$
\|F'(x_*)^{-1}E_F(x,x_*)\|\leq \int_0 ^1 \left \|     F'(x_*)^{-1}[F'(x)-F'(x_*+\tau(x-x_*))]\right\|\,\left\|x_*-x\right\| \; d\tau.
$$
From the last equation and the assumption \eqref{Hyp:MHqn}, we obtain
$$
\|F'(x_*)^{-1}E_F(x,x_*)\| \leq \int_0 ^1  \left[f'\left(\left\|x-x_*\right\|\right)-f'\left(\tau\|x-x_*\|\right)\right]\|x-x_*\|\;d\tau.
$$
Evaluating the above integral and using the definition of $e_f$, the statement follows.
\end{proof}
Define the Newton step to the functions $F$ and $f$ by the following equalities:
\begin{equation} \label{eq:ns}
S_{F}(x):=-F'(x)^{-1}F(x), \qquad s_{f}(t):=-f(t)/f'(t).
\end{equation}
\begin{lemma}  \label{passonewton}
If  $\|x-x_*\|< \min\{\nu, \kappa\}$, then  $
\|S_{F}(x)\|\leq s_{f}(\|x-x_{*}\|). $
\end{lemma}
\begin{proof}
Using \eqref{eq:ns}, $F(x_{*})=0$ and some algebraic manipulation, it follows from   \eqref{eq:def.er} that
\begin{align*}
\|S_{F}(x)\|&=\|-F'(x)^{-1}\left(F(x_{*})-[F(x)+F'(x)(x_{*}-x)]\right)+(x_{*}-x)\|\\
&\leq\|F'(x)^{-1}F'(x_{*})\|\|F'(x_{*})^{-1}\left(F(x_{*})-[F(x)+F'(x)(x_{*}-x)]\right)\|+\|x_{*}-x\|\\
&=\|F'(x)^{-1}F'(x_{*})\|\|F'(x_{*})^{-1}E_{F}(x,x_{*})\|+\|x-x_{*}\|.
\end{align*}
Combining the last equation with Lemma~\ref{wdns} and
Lemma~\ref{pr:taylor} we have 
$$\|S_{F}(x)\|\leq
\frac{e_{f}(\|x-x_{*}\|,0)}{|f'(\|x-x_{*}\|)|}+\|x-x_{*}\|.$$
Since  $f'<0$ in $[0, \,
\nu)$ and $\|x-x_*\|<\nu$, we obtain from last inequality  
\eqref{eq:def.erf}  and  {\bf h1}, that
$$
\|S_{F}(x)\|\leq \frac{f(0)-f(\|x-x_{*}\|)+f'(\|x-x_{*}\|)\|x-x_{*}\|}{-f'(\|x-x_{*}\|)}+\|x-x_{*}\|
=\frac{f(\|x-x_{*}\|)}{f'(\|x-x_{*}\|)}.
$$
So, the last inequality together with the second equality in  \eqref{eq:ns} implies the desired inequality.
 \end{proof}
\begin{lemma} \label{l:wdef}
Let $\Omega\subseteq \mathbb{R}^{n}$ be an open set and $F:{\Omega}\to \mathbb{R}^{n}$ a continuously differentiable function. Let  $x_*\in \Omega$,  $R>~0$ and 
$
\kappa:=\sup \left\{ t\in [0, R): B(x_*, t)\subset\Omega \right\}.
$
Suppose that $F(x_*)=0$, $F '(x_*)$ is invertible and there exists a $f:[0,\; R)\to \mathbb{R}$ continuously differentiable satisfying \eqref{Hyp:MHqn}, {\bf h1} and {\bf h2}.  Let $\vartheta$,  $\omega_{2}$, $\omega_{1}$,  $\nu$, $\rho$ and $\sigma$ as  in Theorem~\ref{th:ntqn}. Assume that   $x\in B(x_*, \sigma)\backslash \{x_*\}$, i.e., $0<\|x-x_*\|< \sigma$.  Define
\begin{equation} \label{eq:DNSqnG}
x_{+}={x}+S, \qquad  B(x)S=-F(x)+r,
\end{equation}
where $B(x)$ is a invertible approximation of $F'(x)$ satisfying 
\begin{equation}\label{con:qnG}
\|B(x)^{-1}F'(x)\| \leq \omega_{1}, \qquad
\|B(x)^{-1}F'(x)-I\| \leq \omega_{2}, 
\end{equation}
and that the  residual  $r$ satisfies 
\begin{equation}\label{eq:ERROqnG}
\|P r\|\leq \theta\|P F(x)\|,
\end{equation}
for some  $\theta$ and $P$   non-negative number and invertible matrix, respectively, satisfying 
\begin{equation}\label{con:vG}
\theta \mbox{cond}(P F'(x))\leq
\vartheta,
\end{equation}
 then $x_{+}$ is well defined and there holds
\begin{equation*}\label{c112G}
\|x_{+}-x_{*}\|\leq
\left[\omega_{1}(1+\vartheta)\frac{|n_{f}(\|x-x_*\|)|}{\|x-x_*\|}+
\omega_{1}\vartheta+\omega_{2} \right]\|x-x_*\|.
\end{equation*}
In particular,
$$
\|x_{+}-x_{*}\|< \|x-x_*\|.
$$
\end{lemma}
\begin{proof}
First note that, as $\|x-x_*\|\leq t<\sigma$, it follows from  Lemma \ref{wdns} that $F'(x)$ is invertible. Now, let $B(x)$ a invertible approximation of it satisfying \eqref{con:qnG}. Thus, $x_{+}$ is well defined. Now, as $F(x_*)=0,$  some simple algebraic manipulation and \eqref{eq:DNSqnG} yield
$$
x_{+}-x_{*}=B(x)^{-1}\big(F(x_{*})-[F(x)+F'(x)(x_{*}-x)]\big)+(B(x)^{-1}F'(x)-I)(x_{*}-x)+B(x)^{-1}{r}.
$$
So, the above equation and  \eqref{eq:def.er} give 
$$
x_{+}-x_{*}=B(x)^{-1}E_{F}(x,x_{*})+(B(x)^{-1}F'(x)-I)(x_{*}-x)+B^{-1}(x){r}.
$$
Again, some algebraic manipulation in the above equation, together with the properties of the norm, imply
\begin{multline*}
\|x_{+}-x_{*}\|\leq\|B(x)^{-1}F'(x)\|\|F'(x)^{-1}F'(x_{*})\|\|F'(x_{*})^{-1}E_{F}(x,x_{*})\|\\
+\|B(x)^{-1}F'(x)-I\|\|x-x_{*}\|+\|B(x)^{-1}F'(x)\|\|F'(x)^{-1}P^{-1}\|\|P r\|.
\end{multline*}
Taking into account the assumptions \eqref{con:qnG} e \eqref{eq:ERROqnG} we obtain from the last equation that
\begin{multline*}
\|x_{+}-x_{*}\|\leq\omega_{1}\|F'(x)^{-1}F'(x_{*})\|\|F'(x_{*})^{-1}E_{F}(x, x_{*})\|+\omega_{2}\|x-x_{*}\|\\
+\omega_{1}\theta \|F'(x)^{-1}P^{-1}\|\|P\,F(x)\|.
\end{multline*}
On the other hand,  \eqref{con:vG} implies $\theta \|(P F'(x))^{-1}\|\|P F'(x)\|\leq
\vartheta$. So,  it is easy to see from \eqref{eq:ns} that
$$
\omega_{1}\theta \|F'(x)^{-1}P^{-1}\|\|P\,F(x)\|\leq
\omega_{1}\theta \|(P'F(x))^{-1}\|\|P F'(x)\|\|S_{F}(x)\|\leq
\omega_{1}\vartheta\|S_{F}(x)\|.
$$
Hence, it follows from the two latter equations that
$$
\|x_{+}-x_{*}\|\leq\omega_{1}\|F'(x)^{-1}F'(x_{*})\|\|F'(x_{*})^{-1}E_{F}(x,x_{*})\|+\omega_{2}\|x-x_{*}\|+\omega_{1}\vartheta\|S_{F}(x)\|.
$$
Combining the last equation with Lemma~\ref{wdns}, Lemma~\ref{pr:taylor} and  Lemma~\ref{passonewton} we conclude that
\[
\|x_{+}-x_{*}\|\leq \omega_{1}\frac{e_{f}(||x-x_{*}||,
0)}{|f'(||x-x_{*}||)|}+\omega_{2}\|x-x_{*}\|+\omega_{1}\vartheta\,s_{f}(||x-x_{*}||).
\]
Now, using  \eqref{eq:def.erf}, \eqref{eq:def.nf},   {\bf h1} and \eqref{eq:ns} we have, by direct calculus, 
$$
\frac{e_{f}(\|x-x_*\|,
0)}{|f'(\|x-x_*\|)|}=|n_{f}(\|x-x_*\|)|, \qquad
s_{f}(||x-x_{*}||)=|n_{f}(\|x-x_*\|)|+\|x-x_*\|.
$$
Therefore, it follows from above inequality and the two latter equalities that 
\[
\|x_{+}-x_{*}\|\leq
\omega_{1}|n_{f}(\|x-x_*\|)|+\omega_{2}\|x-x_*\|+\omega_{1}\vartheta \left(|n_{f}(\|x-x_*\|)|+\|x-x_*\|\right),
\]
which is equivalent to the first inequality of the lemma. 

Because  $x\in B(x_*,\sigma)\backslash \{x_*\}$, i.e., $0<\|x-x_*\|<\sigma$ we obtain  the last inequality of the lemma by combining the first one and Proposition~\ref{pr:incr102} with $t=\|x-x_*\|$.
\end{proof}
\subsection{Proof of  {\bf Teorem \ref{th:ntqn}}} \label{sec:proofqn}
We are now in position to prove Theorem \ref{th:ntqn}.
\begin{proof}
Since $x_0\in B(x_*,\sigma)\backslash \{x_*\}$, i.e., $0<\|x_0-x_*\|<\sigma$, a straighforward induction argument and the last inequality in  Lemma \ref{l:wdef} implies that the sequence $\{x_k\}$ generated by inexact  Newton-like methods is well defined and 
contained in  $B(x_*,\sigma)$.

Our task is now to show that  $\{x_k\}$ converges to $x_*$. Because, $\{x_k\}$ is well defined and  contained in  $B(x_*,\sigma)$, applying Lemma \ref{l:wdef} with  $x_{+}=x_{k+1}$, $x=x_{k}$, $r=r_{k}$, $B(x)=B(x_k)$, $P=P_k$ and $\theta=\theta_k$  we obtain
\begin{equation}\label{c112}
\|x_{k+1}-x_{*}\|\leq
\left[\omega_{1}(1+\vartheta)\frac{|n_{f}(\|x_k-x_*\|)|}{\|x_k-x_*\|}+
\omega_{1}\vartheta+\omega_{2} \right]\|x_k-x_*\|, \qquad \;k=0,1, \ldots.
\end{equation}
$$
\|x_{k+1}-x_*\|<\|x_k-x_*\|, \qquad \;k=0,1, \ldots.
$$
In particular, the last inequality implies that $\|x_{k+1}-x_*\|< \|x_0-x_*\|$, for $k=0,1, \ldots,$ which, together with  \eqref{c112}  and Corollary~\ref{pr:10}, gives
\begin{equation}\label{t1qn}
\|x_{k+1}-x_*\|\leq\left[\omega_{1}(1+\vartheta)\frac{|n_{f}(\|x_{0}-x_*\|)|}{\|x_{0}-x_*\|}+\omega_{1}\vartheta+
\omega_{2} \right] \|x_{k}-x_*\|, \quad  k=0,1,\ldots.
\end{equation}
As Proposition \ref{pr:incr102} with  $t=\|x_0-x_*\|$ gives 
$
\omega_{1}(1+\vartheta)|n_{f}(\|x_{0}-x_*\|)|/\|x_{0}-x_*\|+\omega_{1}\vartheta+\omega_{2}<1,
$
we conclude from above equation that $\{\|x_{k}-x_*\|\}$ converges to zero. So,  $\{x_k\}$ converges to $x_*$.

It remains to prove the last inequality of the theorem. For this, use \eqref{t1qn} and the definition in~\eqref{eq:def.nf}.
\end{proof}
\begin{remark}
If a continuously differentiable function $f:[0,\; \kappa)\to \mathbb{R}$ is a majorant function satisfying the conditions {\bf h1} and  {\bf h2},  then  the function $ h:(-\kappa, \,\kappa)\to \mathbb{R}$ is defined by
\begin{equation} \label{eq:dh}
 h(t)=
      \begin{cases}
       -f(-t), \quad \;\; t\in  (-\kappa, \,0],\\
       \;\;f(t), \quad \quad   \;\; t\in [0, \,\kappa).
      \end{cases}
\end{equation}
 satisfies all hypotheses of Theorem \ref{th:ntqn}. Indeed, it is straightforward to show that $h(0)=0$,  $h'(0)=-1$, $h'(t)=f'(|t|)$ and that
$$
\left|h'(0)^{-1}\left[h'(t)-h'(\tau t)\right]\right| \leq
    f'(|t|)-f'(\tau|t|), \qquad \tau \in [0,1], \quad t\in (-\kappa,\, \kappa).
$$
So, $F=h$, $n=1$,  $\Omega=(-\kappa,\, \kappa)$ and $x_*=0$ satisfy all hypotheses  of Theorem \ref{th:ntqn}. Therefore, we can apply Theorem \ref{th:ntqn} to solve $h(t)=0$.

 Note that if $f'$ is not Lipschitz, then $h'$ is also not Lipschitz. Therefore, we conclude that Theorem~\ref{th:ntqn} enlarges the range of application of theorems on inexact Newton-like Methods having the Lipschitz condition as a hypothesis on the first derivative of the non-linear operator under consideration.
\end{remark} 
Now, we will give some examples of majorant functions satisfying the conditions {\bf h1} and  {\bf h2} with a first derivative that is not Lipschitz.
\begin{example} \label{ex:mf}
 The following functions  satisfy the conditions {\bf h1} and  {\bf h2}:
\begin{itemize}
  \item[{i)}] $f: [0, +\infty)\to \mathbb{R}$ such that $f(t)=e^t-2t-1$;
  \item[{ii)}] $f:[0, \kappa)\to \mathbb{R}$ such that $f(t)=-(1/\kappa)\ln(1-\kappa t)-2t$;
  \item[{iii)}] $f:[0, \kappa)\to \mathbb{R}$ such that $f(t)=[(1-\kappa t)/\kappa]\ln(1-\kappa t)$.
\end{itemize}
Note that the first derivatives of each of the functions above are not Lipschitz. 
\end{example}
\begin{remark}
The assumption \eqref{Hyp:MHqn} was crucial for our analysis. It is worth pointing out that, under appropriate regularity conditions on the nonlinear operator $F$, the assumption \eqref{Hyp:MHqn} always holds in a suitable neighborhood  of $x_*$. For instance, if  $F$ is twice continuously differentiable, then the majorant function $f:[0, \kappa)\to \mathbb{R}$  defined by 
$
f(t)=Kt^{2}/2-t,
$ where $K=\sup\{\|F'(x_*)^{-1}F''(x)\|: x\in B[x_*, \kappa)\}$ satisfies the assumption \eqref{Hyp:MHqn}.  Estimating the constant $K$ is a very difficult problem. Therefore, the goal is to identify  classes of nonlinear operators for which it is possible to obtain a majorant function. We will give some examples of such  classes in the next section.
\end{remark} 
\section{Special cases} \label{sec:ec}
In this section we present three special cases of Theorem~\ref{th:ntqn}. Namely, convergence results under an affine invariant Lipschitz condition, Smale's condition for
analytical functions and   Nesterov-Nemirovskii's condition  for  self-concordant functions.

\subsection{Convergence result for affine invariant Lipschitz condition}
In this section we show a correspondent theorem to Theorem \ref{th:ntqn} under an affine invariant Lipschitz condition (see \cite{DH}, \cite{C08}  and  \cite{B10}) instead of the general assumption \eqref{Hyp:MHqn}.

\begin{theorem}\label{th:ntqnnm}
Let $\Omega\subseteq \mathbb{R}^{n}$ be an open set and
$F:{\Omega}\to \mathbb{R}^{n}$ be continuously differentiable in
$\Omega$. Take  $x_*\in \Omega$ and let
$$
\kappa:=\sup \left\{ t>0: B(x_*, t)\subset\Omega \right\}.
$$
Assume that $F '(x_*)$ is invertible, $F(x_*)=0$, and there exists a
$K>0$  such that
\begin{equation}\label{eq:ailc}
\left\|F'(x_*)^{-1}\left[F'(x)-F'(y)\right]\right\| \leq K\|x-y\|, \qquad \forall\; x, y\in B(x_*, \kappa).
\end{equation}
Take $0\leq \vartheta<1$, $0\leq \omega_{2}<\omega_{1}$  such that  $\omega_{1}\vartheta+\omega_{2}< 1$.  Let
$$
\sigma:=\min\left\{\kappa,
\,2(1-\vartheta\omega_{1}-\omega_{2})/\left(K(2+\omega_1-\vartheta\omega_1-2\omega_2)\right)
\right\}.
$$
 Then, the inexact Newton-like method for solving
$F(x)=0$, with an initial point $x_0\in B(x_*, \sigma) \backslash \{x_*\}$
\begin{equation}\label{eq:ailcni}
x_{k+1}={x_k}+S_k, \qquad  B(x_k)S_k=-F(x_k)+r_{k}, \qquad
\; k=0,1,\ldots,
\end{equation}
where $B(x_k)$ is an invertible approximation of $F'(x_k)$
satisfying
$$
\|B(x_k)^{-1}F'(x_k)\| \leq \omega_{1}, \qquad
\|B(x_k)^{-1}F'(x_k)-I\| \leq \omega_{2},
$$
and the  residual  $r_k$ satisfies
\begin{equation}\label{eq:ailrc}
\|P_{k}r_{k}\|\leq \theta_{k}\|P_{k}F(x_{k})\|,
\end{equation}
for some invertible matrix sequence $\{P_{k}\}$ of preconditioners and 
forcing sequence $\{\theta_{k}\}$   of non-negatives numbers, satisfying 
$$
\theta_{k}\mbox{cond}(P_{k}F'(x_{k}))\leq \vartheta,
$$
is well defined, contained in $B(x_*,\sigma)$, converges to $x_*$
and there holds
$$
    \|x_{k+1}-x_*\| \leq
\left[\omega_1(1+\vartheta)\frac{K\|x_{0}-x_*\|}{2(1-K\|x_{0}-x_*\|)}+\omega_1\vartheta+\omega_2
\right]\|x_k-x_*\|, \quad \forall \; k=0,1,\ldots .
 $$
 \end{theorem}
 \begin{proof}
It is immediately possible to prove that  $F$, $x_*$ and $f:[0, \kappa)\to
\mathbb{R}$ defined by $ f(t)=Kt^{2}/2-t, $ satisfy the inequality
\eqref{Hyp:MHqn} and the conditions  {\bf h1} and  {\bf h2} in
Theorem \ref{th:ntqn}. In this case, it is easy to see that the
constants $\rho$ and $\nu$, as defined in Theorem \ref{th:ntqn},
satisfy
$$\rho=2(1-\vartheta\omega_{1}-\omega_{2})/\left(K(2+\omega_1-\vartheta\omega_1-2\omega_2)\right) \leq \nu=1/K,$$
as a consequence
$$
\sigma:=\min \{\kappa,\;
2(1-\vartheta\omega_{1}-\omega_{2})/\left(K(2+\omega_1-\vartheta\omega_1-2\omega_2)\right)\}.
$$
Therefore, as  $F$, $\sigma$,  $f$ and $x_*$ satisfy all of the
hypotheses of  Theorem \ref{th:ntqn}, taking  $x_0\in B(x_*,
\sigma)\backslash \{x_*\}$ the statements of the theorem follow from
Theorem \ref{th:ntqn}.
\end{proof}
Although the condition  \eqref{eq:ailc} is affine invariant (it is insensitive with respect to transformation  of the map $F$ of the form $F\mapsto AF$), iteration \eqref{eq:ailcni} and the condition for the residual \eqref{eq:ailrc} is not affine invariant. So,  Theorem~\ref{th:ntqnnm} is not affine invariant. Now, taking $B(x_k)=F'(x_k)$ in iteration \eqref{eq:ailcni} and $P_k=F'(x_k)^{-1}$ in the condition for the residual \eqref{eq:ailrc}, the Theorem~\ref{th:ntqnnm} becomes  affine invariant. It is easy to see that, for the theorem that uses the Lipschitz condition
$$
\left\|F'(x)-F'(y)\right\| \leq L\|x-y\|, \qquad \forall\; x, y\in B(x_*, \kappa),
$$
instead of the affine invariant Lipschitz condition \eqref{eq:ailc}, the convergence radius is given by
$$
\sigma:=\min\left\{\kappa,
\,2(1-\vartheta\omega_{1}-\omega_{2})/\left(L\|F'(x_*)\|(2+\omega_1-\vartheta\omega_1-2\omega_2)\right)
\right\}.
$$
We point out that the convergence radius of affine invariant  theorems are
 insensitive to invertible transformation of the map $F$, but  that theorems with the Lipschitz condition (see next example) are sensitive. For more details about  affine invariant theorem  see \cite{DH}.
\begin{example} Assume that $B(x_k)=F'(x_k)$,  $P_k=F'(x_k)^{-1}$, $\omega_1=1$ and $\omega_2=0$ in  Theorem \ref{th:ntqnnm}.
Let $F: {\mathbb{R}} {^2}\to {\mathbb{R}} {^2}$ be given by $F(x_1, x_2)=(x_1^2/2-x_1, \; x_2^2/2-x_2)$. Note that $ F(0, 0)= (0, 0)$.
Using the Euclidean vector norm and the associated operators norm, it is easy to see that $\|F''(x_1, x_2)\|=1$ and the Lipschitz constant for $F'$ is $1$.   In this case, the Lipschitz condition for $F'$ and the affine invariant Lipschitz condition \eqref{eq:ailc}  are equal. Therefore, by applying Theorem \ref{th:ntqnnm} we conclude that the convergence radius for solving $F(x_1, x_2)=0$ is $2(1-\vartheta)/\left(3-\vartheta\right).$
Let the invertible matrix
$$
A=\left[      
\begin{matrix}
1&0 \\
0&1/\epsilon
\end{matrix}
\right], \qquad  0<\epsilon <1,
$$
and the map  $G:\Omega\to {\mathbb{R}} {^n}$ given by $G(x)=AF(x).$ Hence $G'(x_1, x_2)=AF'(x_1, x_2)$. Moreover, 
$$
\|G'(0, 0)^{-1}\|=1, \qquad \|G''(x_1, x_2)\|=1/\epsilon, 
$$
and the Lipschitz constant for $G'$ is $1/\epsilon$. Applying Theorem~\ref{th:ntqnnm} with $B(x_k)=G'(x_k)$ and  $P_k=G'(x_k)^{-1}$, we conclude that,  due to its insensitivity  to invertible transformation,  the convergence radius for solving $G(x_1, x_2)=0$ is also $2(1-\vartheta)/\left(3-\vartheta\right)$. However, if in Theorem~\ref{th:ntqnnm} the Lipschitz condition for $F'$ is assumed instead of the affine invariant Lipschitz condition \eqref{eq:ailc}, then the convergence radius   is $(2\epsilon)(1-\vartheta)/\left(3-\vartheta\right)$.
\end{example}
\subsection{Convergence result under Smale's condition }
In this section we show a  correspondent theorem to Theorem
\ref{th:ntqn} under  Smale's condition. For more details about
 Smale's condition   see \cite{S86}.

\begin{theorem}\label{theo:Smale}
Let $\Omega\subseteq \mathbb{R}^{n}$ be an open set and
$F:{\Omega}\to \mathbb{R}^{n}$ an analytic function. Take $x_*\in
\Omega$ such that $F '(x_*)$  is invertible and  $F(x_*)=0$  and let

\begin{equation} \label{eq:SmaleCond}
\kappa:=\sup\{t >0 : B(x_*, t)\subset \Omega\} \qquad
\mbox{and}\qquad \gamma := \sup _{ n > 1 }\left\| \frac
{F'(x_*)^{-1}F^{(n)}(x_*)}{n !}\right\|^{1/(n-1)}<+\infty.
\end{equation}
Take $0\leq \vartheta<1$, $0\leq \omega_{2}<\omega_{1}$  such that  $\omega_{1}\vartheta+\omega_{2}< 1$.  Let $a=\omega_1(1+\vartheta)$, 
$b=(1-\omega_1\vartheta-\omega_2)$ and
$$
\sigma:=\min \left\{\kappa,\frac{ a+4b-( \sqrt{(a+4b)^2-8b^2})}{4b\gamma}\right\}.
$$
Then, the inexact Newton-like method for solving
$F(x)=0$, with initial point $x_0\in B(x_*, \sigma) \backslash \{x_*\}$
$$
x_{k+1}={x_k}+S_k, \qquad  B(x_k)S_k=-F(x_k)+r_{k}, \qquad
\; k=0,1,\ldots,
$$
where $B(x_k)$ is an invertible approximation of $F'(x_k)$
satisfying
$$
\|B(x_k)^{-1}F'(x_k)\| \leq \omega_{1}, \qquad
\|B(x_k)^{-1}F'(x_k)-I\| \leq \omega_{2},
$$
and the  residual  $r_k$ satisfies
$$
\|P_{k}r_{k}\|\leq \theta_{k}\|P_{k}F(x_{k})\|,
$$
for some forcing sequence $\{\theta_{k}\}$   of non-negative numbers and
an invertible matrix sequence $\{P_{k}\}$ of preconditioners, satisfying 
$$
\theta_{k}\mbox{cond}(P_{k}F'(x_{k}))\leq \vartheta,
$$
is well defined, contained in $B(x_*,\sigma)$, and converges to $x_*$ 
$$
    \|x_{k+1}-x_*\| \leq
\left[\omega_1(1+\vartheta)\frac{\gamma\|x_{0}-x_*\|}{2(1-\gamma\|x_{0}-x_*\|)^2-1}+\omega_1\vartheta+\omega_2
\right]\|x_k-x_*\|, \quad \forall \; k=0,1,\ldots .
 $$
 \end{theorem}
We need the following result to prove the above theorem.
\begin{lemma} \label{lemma:qc1}
Let $\Omega\subseteq \mathbb{R}^{n}$ be an open set and
$F:{\Omega}\to \mathbb{R}^{n}$ an analytic function.  Suppose that
$x_*\in \Omega$,  $F '(x_*)$ is  invertible and that $B(x_{*},
1/\gamma) \subset \Omega$, where $\gamma$ is defined in
\eqref{eq:SmaleCond}. Then, for all $x\in B(x_{*}, 1/\gamma)$
there holds
$$
\|F'(x_*)^{-1}F''(x))\| \leqslant  (2\gamma)/( 1- \gamma
\|x-x_*\|)^3.
$$
\end{lemma}
\begin{proof}
Let $x\in \Omega$.  Since $F$ is an analytic function, we have
$$
F'(x_{*})^{-1}F''(x)= \sum _{n=0}^{\infty}\frac
{1}{n!}F'(x_{*})^{-1}F^{(n + 2)}(x_{*})(x - x_{*})^{n}.
$$
Combining  \eqref{eq:SmaleCond} and  the above equation we obtain,
after some simple calculus, that
$$
\|F'(x_{*})^{-1}F''(x)\| \leqslant \,\gamma \sum
_{n=0}^{\infty}(n+2)(n+1)(\gamma ||x-x_{*}||)^{n} .
$$
On the other hand, as $B(x_{*}, 1/\gamma) \subset \Omega$ we have
 $\gamma \|x-x_*\|< 1$. So,  from Proposition~\ref{le:ess}  we
conclude
$$
\frac{2}{(1-\gamma\|x-x_*\|)^3}=\sum_{n=0}^{\infty}(n+2)(n+1)(\gamma
||x-x_{*}||)^{n}.
$$
Combining the two above equations, we obtain the desired
result.
\end{proof}
The next result gives a condition that is easier to check than
condition \eqref{Hyp:MHqn}, whenever the functions under consideration
are twice continuously differentiable.
\begin{lemma} \label{lc}
Let $\Omega\subseteq \mathbb{R}^{n}$ be an open set and
  $F:{\Omega}\to \banachb$ be twice continuously on $\Omega$.
 Let $x_*\in \Omega$ with $F '(x_*)$  be invertible.  If there exists a \mbox{$f:[0,R)\to \mathbb {R}$} twice continuously differentiable such that
 \begin{equation} \label{eq:lc2}
\|F'(x_*)^{-1}F''(x)\|\leqslant f''(\|x-x_*\|),
\end{equation}
for all $x\in  \Omega$ such that  $\|x-x_*\|<R$. Then $F$ and $f$
satisfy \eqref{Hyp:MHqn}.
\end{lemma}
\begin{proof}
   Taking $\tau \in [0,1]$ and $x\in \Omega$, such that $x_*+\tau(x-x_*)\in \Omega$  and   $\|x-x_*\|<R$, we
  obtain that
 $$
 \|F'(x_*)^{-1}\left[F'(x)-F'(x_*+\tau(x-x_*))\right]\|\leq
  \int_{\tau}^{1}\|F'(x_*)^{-1}F''(x_*+t(x-x_*))\|\,\|x-x_*\|dt.
$$
  Now, as $\|x-x_*\|<R$ and $f$ satisfies \eqref{eq:lc2}, we obtain from the last inequality that
\begin{align*}
 \|F'(x_*)^{-1}\left[F'(x)-F'(x_*+\tau(x-x_*))\right]\|&\leq \int_{\tau}^{1}f''(t\|x-x_*\|)\|x-x_*\|dt.
\end{align*}
Evaluating the latter integral, the statement follows.
\end{proof}
{\bf [Proof of Theorem \ref{theo:Smale}]}. Assume that all
hypotheses of   Theorem \ref{theo:Smale} hold. Consider the
real function $f:[0,1/\gamma) \to \mathbb{R}$ defined by
$$
f(t)=\frac{t}{1-\gamma t}-2t.
$$
It is straightforward to show that $f$ is  analytic and that
$$
f(0)=0, \quad f'(t)=1/(1-\gamma t)^2-2, \quad f'(0)=-1, \quad
f''(t)=(2\gamma)/(1-\gamma t)^3, \quad f^{n}(0)=n!\,\gamma^{n-1},
$$
for $n\geq 2$. From  the last four equalities it is easy to see that
$f$ satisfies {\bf h1}  and  {\bf h2}. Now, since $f''(t)=(2\gamma)/(1-\gamma t)^3$ combining
Lemma~\ref{lc}, Lemma~\ref{lemma:qc1}  we
conclude that $F$  and $f$ satisfy  \eqref{Hyp:MHqn} with $R=1/\gamma$. Define
$$
\nu:=\sup\{t\in [0, 1/\gamma): f'(t)<0\}, \quad \rho:=\sup\{t\in(0, \nu):\omega_{1}(1+\vartheta)[f(t)/(tf'(t))-1]+\omega_{1}\vartheta+\omega_{2}<1\}.
$$
In this case, it is easy to see that the constants $\nu$ and $\rho$  satisfy
$$\rho=\frac{ a+4b-( \sqrt{(a+4b)^2-8b^2})}{4b}, \qquad \nu=\frac{\sqrt{2}-1}{\sqrt{2}\gamma}, \qquad \rho<\nu< \frac{1}{\gamma}, 
$$
where $a=\omega_1(1+\vartheta)$, 
$b=(1-\omega_1\vartheta-\omega_2)$. Finally, let $\sigma:=\min\{\kappa,\rho\}.$ Therefore, as  $F$, $\sigma$, 
$f$ and $x_*$ satisfy all hypothesis of  Theorem \ref{th:ntqn},
taking $x_0\in B(x_*, \sigma)\backslash \{x_*\}$,  the statements of the
theorem follow from Theorem \ref{th:ntqn}.  \qed
\subsection{Convergence result under The Nesterov-Nemirovskii condition}
In this section we show a correspondent theorem to Theorem
\ref{th:ntqn} under  the  Nesterov-Nemirovskii condition(see \cite{NN-94}). 

Let $\Omega\subset \mathbb{R}^n$ be a convex set. A  function
$g:\Omega \to \mathbb{R}$ is called $a$-self-concordant with the
parameter $a>0$,  if $g\in C^{3}(\Omega)$ , i.e.,  three times
continuously differentiable in $\Omega$,  is a  convex function on
$\Omega$ and satisfies the following inequality
\begin{equation} \label{eq:scf}
|g'''(x)[h, h, h]|\leqslant 2a^{-1/2}(g''(x)[h, h])^{3/2}, \qquad
\forall \;x\in \Omega,\; \forall\;h\in \mathbb{R}^n.
\end{equation}
Take $x_*\in \Omega$ such that  $g''(x_*)$ is invertible. Define
$\banacha:=(\mathbb{R}^n, \langle ., .
\rangle_{x_*})$ as  the Euclidean  space $\mathbb{R}^n$ with the
inner product  and the associated norm defined, respectively,  by
\begin{equation*}
\langle u, v \rangle_{x_*}:=a^{-1}\langle g''(x_*) u, v \rangle,
\qquad \|u\|_{x_*}:=\sqrt{\langle u, u \rangle_{x_*}}, 
\end{equation*}
where $\langle ., .\rangle$ is the Euclidean inner product.
So,  the open and closed ball of radius $r>0$  centered at $x_*$ (
Dikin's ellipsoid of radius $r>0$  centered at $x_*$ ) in
$\banacha$ are defined, respectively,  as
$$
W_{r}(x_*):=\left\{x\in \mathbb{R}^n:
\|x-x_*\|_{x_*}<r\right\},\qquad W_{r}[x_*]:=\left\{x\in
\mathbb{R}^n: \|x-x_*\|_{x_*}\leq r \right\}.
$$
\begin{theorem}\label{theo:Self-Conc}
Let $\Omega\subseteq \banacha$ be a convex  set and
$g:\Omega \to \mathbb{R}$  an $a$-self-concordant function.
Take $x_*\in \Omega$ with $g''(x_*)$  invertible and let $\kappa:=\sup\{u >0 : W_{u}(x_*)\subset \Omega\}$.
  Suppose that $g'(x_*)=0$.\\
   Take $0\leq \vartheta<1$, $0\leq \omega_{2}<\omega_{1}$  such that  $\omega_{1}\vartheta+\omega_{2}< 1$.  Let $a=\omega_1(1+\vartheta)$, 
$b=(1-\omega_1\vartheta-\omega_2)$ and
$$
\sigma:=\min \left\{\kappa,\frac{ a+4b-( \sqrt{(a+4b)^2-8b^2})}{4b}\right\}.
$$
Then, the inexact Newton-like method for solving
$g'(x)=0$, with an initial point $x_0\in B(x_*, \sigma)\backslash \{x_*\}$
$$
x_{k+1}={x_k}+S_k, \qquad  B(x_k)S_k=-F(x_k)+r_{k}, \qquad
\; k=0,1,\ldots,
$$
where $B(x_k)$ is an invertible approximation of $g''(x_k)$
satisfying
$$
\|B(x_k)^{-1}g''(x_k)\| \leq \omega_{1}, \qquad
\|B(x_k)^{-1}g''(x_k)-I\| \leq \omega_{2},
$$
and the  residual  $r_k$ satisfies
$$
\|P_{k}r_{k}\|\leq \theta_{k}\|P_{k}F(x_{k})\|,
$$
for some forcing sequence $\{\theta_{k}\}$   of non-negative numbers and
an invertible matrix sequence $\{P_{k}\}$ of preconditioners, satisfying 
$$
\theta_{k}\mbox{cond}(P_{k}F'(x_{k}))\leq \vartheta,
$$
is well defined, contained in $B(x_*,\sigma)$, converges to $x_*$
and $$
    \|x_{k+1}-x_*\| \leq
\left[\omega_1(1+\vartheta)\frac{\|x_{0}-x_*\|}{2(1-\|x_{0}-x_*\|)^2-1}+\omega_1\vartheta+\omega_2
\right]\|x_k-x_*\|, \quad \forall \; k=0,1,\ldots .
 $$ \end{theorem}
We need some auxiliary results about self-concordant functions to
prove the above theorem. We begin with two well known propositions
in the theory of self-concordant functions,  from  Nesterov and
Nemirovskii \cite{NN-94}.
\begin{proposition} \label{scf2}
Let $\Omega\subset \banacha$ be an open convex set and  let
$g:\Omega \to \mathbb{R}$ be an $a$-self-concordant function.
Then,
$$
|g'''(x)[h_1, h_2, h_3]|\leqslant
2a^{-1/2}\Pi_{i=1}^{3}(g''(x)[h_{i}, h_{i}])^{1/2},\qquad
\forall\; x\in \Omega,\; \forall\; h_1, h_2, h_3\in \banacha.
$$
\end{proposition}
\begin{proof}
See Proposition 9.1.1, Appendix 1,  pp.361 of \cite{NN-94}.
\end{proof}
\begin{proposition} \label{scf3}
Let $\Omega\subset \banacha$ be an open convex set and  let
$g:\Omega \to \mathbb{R}$ be an $a$-self-concordant function.
Assume that $W_{1}(x_*)\subset \Omega$. Then there holds
\begin{equation*}
g''(x)[h, h]\leq \frac{1}{(1-\|x-x_*\|_{x_*})^{2}}g''(x_*)[h, h],
\qquad \forall \;x\in W_{1}(x_*), \quad \forall \;h\in
\banacha.
\end{equation*}
\end{proposition}
\begin{proof}
See Theorem 2.1.1 pp.13 of \cite{NN-94}.
\end{proof}
The next result is a combination of the two last propositions, which has
appeared  in  \cite{ABM2004} Lemma~5.1. We include
the proof here.
\begin{lemma} \label{lipself}
Let $\Omega\subset \banacha$ be an open convex set and  let
$g:\Omega \to \mathbb{R}$ be an  $a$-self-concordant function.
Assume that $W_{1}(x_*)\subset \Omega$. Then  
\begin{equation*}
\|g''(x_*)^{-1}g'''(x)\|_{x_*}\leq
\frac{2}{(1-\|x-x_*\|_{x_*})^{3}}, \qquad \forall \;x\in
W_{1}(x_*).
\end{equation*}
\end{lemma}
\begin{proof}
Letting  $x\in W_{1}(x_*)$ and $h_{1}$, $h_{2}$, $h_{3}\in
\mathbb{R}^n$ we have from \eqref{eq:scf} and Proposition
\ref{scf2}  
\begin{align*}
|\langle g''(x_*)^{-1}g'''(x)h_1 h_2, h_3\rangle_{x_*}|
&=  a^{-1}|\langle g''(x_*)\left(g''(x_*)^{-1}g'''(x)\right)h_1 h_2, h_3\rangle | \\
&=a^{-1}|g'''(x)[h_{1}, h_{2}, h_{3}]| \\
&\leq 2 a^{-3/2} \Pi_{i=1}^{3}(g''(x)[h_{i}, h_{i}])^{1/2}.
\end{align*}
Since $\|g''(x_*)^{-1}g'''(x)\|_{x_*}:=\sup \left\{|\langle
g''(x_*)^{-1}g'''(x)h_1 h_2, h_3\rangle_{x_*}|:
\|h_i\|_{x_*}\leqslant 1, \; i=1,2,3 \right\}$, we have from the last
inequality
\begin{equation} \label{eq:lipself1}
\|g''(x_*)^{-1}g'''(x)\|_{x_*}\leq 2 a^{-3/2}\sup
\left\{\Pi_{i=1}^{3}(g''(x)[h_{i}, h_{i}])^{1/2}:
\|h_i\|_{x_*}\leq 1, \; i=1,2,3  \right\}.
\end{equation}
Therefore, it follows from \eqref{eq:lipself1} and  Proposition
\ref{scf3} that
\begin{align*}
 \|g''(x_*)^{-1}g'''(x)\|_{x_*} &\leq \frac{2 a^{-3/2}}{(1-\|x-x_*\|_{x_*})^{3}}\sup \left\{\Pi_{i=1}^{3}g''(x_*)[h_{i}, h_{i}])^{1/2}: \|h_i\|_{x_*}\leq 1, \; i=1,2,3  \right\} \notag\\
                         &= \frac{2 }{(1-\|x-x_*\|_{x_*})^{3}}\sup \left\{\Pi_{i=1}^{3}a^{-1/2}g''(x_*)[h_{i}, h_{i}])^{1/2}: \|h_i\|_{x_*}\leq 1, \; i=1,2,3  \right\} \notag\\
                         &=\frac{2 }{(1-\|x-x_*\|_{x_*})^{3}}\sup \left\{\Pi_{i=1}^{3}\|h_i\|_{x_*} : \|h_i\|_{x_*}\leq 1, \; i=1,2,3  \right\}\\
                        &\leq \frac{2 }{(1-\|x-x_*\|_{x_*})^{3}},
\end{align*}
which is a proof of the Lemma.
\end{proof}
{\bf [Proof of Theorem \ref{theo:Self-Conc}]}.  Assume that all
hypotheses of  Theorem \ref{theo:Smale} hold. Consider the
real function $f:[0, 1) \to \mathbb{R}$ defined by
$$
f(t)=\frac{t}{1- t}-2t.
$$
It is straightforward to show that $f$ is  analytic and that
$$
f(0)=0, \quad f'(t)=1/(1- t)^2-2, \quad f'(0)=-1, \quad
f''(t)=2/(1-t)^3, \quad f^{n}(0)=n!,
$$
for $n\geq 2$. From the last four equalities, it is easy to conclude that
$f$ satisfies {\bf h1}  and  {\bf h2}. Now, combining
Lemma~\ref{lc}, Lemma~\ref{lipself} and the latter equality we
obtain that $g'$  and $f$ satisfy  \eqref{Hyp:MHqn} with $R=1$. Define
$$
\nu:=\sup\{t\in [0, 1): f'(t)<0\}, \quad \rho:=\sup\{t\in(0, \nu):\omega_{1}(1+\vartheta)[f(t)/(tf'(t))-1]+\omega_{1}\vartheta+\omega_{2}<1\}.
$$
In this case, it is easy to see that the constants $\nu$ and $\rho$  satisfy
$$\rho=\frac{ a+4b-( \sqrt{(a+4b)^2-8b^2})}{4b}, \qquad \nu=\frac{\sqrt{2}-1}{\sqrt{2}}, \qquad \rho<\nu< 1, 
$$
where $a=\omega_1(1+\vartheta)$, 
$b=(1-\omega_1\vartheta-\omega_2)$. Finally, let $\sigma:=\min\{\kappa,\rho\}.$ Therefore, as  $F=g'$, $\sigma$, 
$f$ and $x_*$ since the above satisfy all hypotheses of  Theorem \ref{th:ntqn},
taking $x_0\in B(x_*, \sigma)\backslash \{x_*\}$,  the statements of the
theorem follow from Theorem \ref{th:ntqn}.  \qed

\section{Final remarks} \label{rf}

As pointed out by  Morini in  \cite{B10} if preconditioning  $P_{k}$,  satisfying  
\begin{equation} \label{eq:r1}
\|P_{k}r_{k}\|\leq \theta_{k}\|P_{k}F(x_{k})\|,
\end{equation}
for some forcing sequence $\{\theta_{k}\}$, is applied  in finding the  inexact Newton steep, then the inverse proportionality  between each forcing term $\theta_{k}$ and $\mbox{cond}(P_{k}F'(x_{k}))$ stated in the following assumption: 
\begin{equation} \label{eq:r2}
0< \theta_{k}\mbox{cond}(P_{k}F'(x_{k}))\leq
\vartheta,\qquad 
\; k=0,1,\ldots,
\end{equation}
is sufficient to guarantee convergence, and may be overly restrictive to bound the sequence $\{\theta_{k}\}$,  always such that the matrices $P_{k}F'(x_{k})$, for $k=0,1,\ldots,$ are badly conditioned. Moreover,   $\theta_{k}$ does not depend on $\mbox{cond} (P_{k})$ but only on the $\mbox{cond}(P_{k}F'(x_{k}))$ and a suitable choice of scaling matrix $P_k$ leads to a relaxation of the forcing terms. 

Using the assumptions \eqref{eq:r1} and \eqref{eq:r2},  we  presented a  new local convergence analysis for inexact Newton-like methods  under  majorant condition.  In our analysis, the  affine invariant Lipschitz condition (see \cite{DH}, \cite{C08} and  \cite{B10})  is relaxed  by using the  majorant condition (see equation  \eqref{Hyp:MHqn}  in Theorem \ref{th:ntqn}). Although the condition \eqref{Hyp:MHqn} is equivalent to the Chen and Li condition (see  equation $(1.4)$  in \cite{C10}),  our analysis elucidates the relationship of the majorant function  with the non-linear operator under consideration (see Lemma~\ref{wdns}. In addition, Lemma~\ref{pr:taylor},  Lemma~\ref{passonewton} and  Lemma~\ref{l:wdef})  allow  us to obtain the  special cases  Theorem~\ref{th:ntqnnm}, Theorem~\ref{theo:Smale} and Theorem~\ref{theo:Self-Conc} of  Theorem \ref{th:ntqn} as an application.

Finally, we point out that the Kantorovich analysis produced a semilocal convergence result, in that it ensures convergence of Newton's Method under very mild assumptions and proves the existence of a solution. On the other hand, local analysis gives us the optimal convergence radius.


\end{document}